\documentstyle[12pt]{article}
\font\tenof=msym10 at 12pt
\def\R{\mbox{\tenof R}}

\def\C{\mbox{\tenof C}}
\def\N{\mbox{\tenof N}}
\def\cR{\mbox{${\cal R}$}}
\def\id{\mbox{\rm id}}
\def\su{\mbox{\rm su$_q$(2)}}
\def\alg{\mbox{${\cal A}^+_q(1)$}}
\def\case#1#2{{\textstyle{#1\over #2}}}
\def\ids{\mbox{\rm\scriptsize id}}
\def\la{\longleftarrow}
\def\ra{\longrightarrow}
\def\da{\Big\downarrow}
\def\ua{\Big\uparrow}
\def\lra{\longleftrightarrow}
\def\uda{\Big\updownarrow}
\def\ss{\scriptstyle}
\def\A{\mbox{${\cal A}$}}
\def\H{\mbox{${\cal H}$}}

\renewcommand{\theequation}{\arabic{section}.\arabic{equation}}

\begin{document}
\baselineskip=22pt plus 1pt minus 1pt
%
%
\begin{center}
{\Large \bf
A nonlinear deformed su(2) algebra with a two-colour quasitriangular Hopf
structure}

\bigskip 

{D.~Bonatsos\footnote{ \rm E-mail: bonat at cyclades.nrcps.ariadne-t.gr}}\\
{\it Institute of Nuclear Physics, NCSR Demokritos, GR-15310
Aghia Paraskevi, Attiki, Greece}
\bigskip

{C.~Daskaloyannis\footnote{ \rm E-mail: daskaloyanni@olymp.ccf.auth.gr}}\\
{\it Department of Physics, Aristotle University of
Thessaloniki}\\ 
{\it GR-54006 Thessaloniki, Greece}
\bigskip

{P.~Kolokotronis}\\
{\it Institute of Nuclear Physics, NCSR Demokritos, GR-15310
Aghia Paraskevi, Attiki, Greece}
\bigskip

{A.~Ludu\footnote{ \rm E-mail: ludu@roifa.ifa.ro}}\\
{\it Department of Theoretical Physics, Faculty of Physics,
University of Bucharest}\\ 
{\it Bucharest-Magurele, P.O.~Box MG-5211, Romania}
\bigskip

{C.~Quesne\footnote{ \rm Directeur de recherches FNRS; E-mail:
cquesne@ulb.ac.be}}\\
{\it Physique Nucl\'eaire Th\'eorique et Physique Math\'ematique, 
Universit\'e Libre de Bruxelles}\\
{\it Campus de la Plaine CP229, Boulevard~du Triomphe, B-1050 Brussels,
Belgium}\\ 

\end{center}
\newpage
%
%
\noindent{\bf Abstract.}
Nonlinear deformations of the enveloping algebra of su(2), involving two arbitrary
functions of $J_0$ and generalizing the Witten algebra, were introduced some time
ago by Delbecq and Quesne. In the present paper, the problem of endowing some of
them with a Hopf algebraic structure is addressed by studying in detail a specific
example, referred to as \alg. This algebra is shown to possess two series of
($N+1$)-dimensional unitary irreducible representations, where
$N=0$, 1, 2,~$\ldots$. To allow the coupling of any two such representations, a
generalization of the standard Hopf axioms is proposed by proceeding in two steps.
In the first one, a variant and extension of the deforming functional technique is
introduced: variant because a map between two deformed algebras, \su\ and \alg, is
considered instead of a map between a Lie algebra and a deformed one, and
extension because use is made of a two-valued functional, whose inverse is
singular. As a result, the Hopf structure of \su\ is carried over to \alg, thereby
endowing the latter with a double Hopf structure. In the second step, the definition
of the coproduct, counit, antipode, and \cR-matrix is extended so that the double
Hopf algebra is enlarged into a new algebraic structure. The latter is referred to as
a two-colour quasitriangular Hopf algebra because the corresponding \cR-matrix is
a solution of the coloured Yang-Baxter equation, where the `colour' parameters take
two discrete values associated with the two series of finite-dimensional
representations.\par 
%
%
\vspace{3cm}
\noindent
Running title: Two-colour quasitriangular Hopf algebra

\par
\newpage
%
%
\section{INTRODUCTION} \label{sec:intro}
Quantized universal enveloping algebras, also called $q$-algebras, refer to
some specific deformations of (the universal enveloping algebra of) Lie
algebras, to which they reduce when the deformation parameter~$q$ 
is set equal to one \cite{Drinfeld}. 
The simplest example of $q$-algebra, \su\ $\equiv$
U$_q$(su(2)), was first introduced by Sklyanin~\cite{Sklyanin},
and independently by Kulish and Reshetikhin~\cite{Kulish} 
in their work on Yang-Baxter equations. A
Jordan-Schwinger realization of \su\ in terms of $q$-bosonic operators was
then derived by Biedenharn~\cite{Biedenharn} and 
Macfarlane~\cite{Macfarlane}. Since then, \su\ has
been applied in various branches of physics. It has been found suitable, for
instance, for the solution of deformed spin-chain models~\cite{Pasquier},
as well as for
the approximate description of rotational spectra of deformed
nuclei~\cite{Bon1},
superdeformed nuclei~\cite{Bon2}, and diatomic molecules~\cite{Bon3}.\par
%
%
In addition to the usual version of the deformed su(2) algebra, namely
\su, several generalized forms of the algebra have been introduced.
Deformations involving one arbitrary function of $J_0$ were independently
proposed by Polychronakos~\cite{Poly} and Ro\v cek~\cite{Rocek}.
Their representation theory
is characterized by a rich variety of phenomena, which might be of interest in
applications to particle physics. Bonatsos~{\it et al\/} \cite{Bon4}
considered a
subclass of these algebras having a representation theory as close as possible
to the usual su(2), so that they might prove useful in applications to nuclear
and molecular physics similar to those above-mentioned in connection with
\su.\par
%
%
Deformations of su(2) involving two arbitrary functions of $J_0$ were
introduced by Delbecq and Quesne~\cite{Quesne}.
Contrary to the former deformations, for
which the spectrum of $J_0$ is linear as for su(2), the latter give rise to
exponential spectra. Such spectra did recently arouse much interest in various
contexts, for instance in connection with alternative Hamiltonian
quantization~\cite{Fairlie},
exactly solvable potentials~\cite{Spiri1}, $q$-deformed
supersymmetric quantum mechanics~\cite{Spiri2}, and q-deformed interacting
boson models~\cite{Gupta}.\par
%
%
From a mathematical viewpoint, $q$-algebras are just special classes of Hopf
algebras \cite{Drinfeld,Majid}.
One of the basic data defining a Hopf algebra is the so-called
comultiplication rule. The latter plays an important role in representation
theory since it allows one to define a product of two independent
representations that is still a representation. Hence, for \su, for
instance, Wigner-Racah calculus can be developed in much the same way as for
the standard su(2) Lie algebra (see e.g.\ Ref.~\cite{Nomura}
and references quoted therein).\par
%
%
The existence of a comultiplication rule for more general deformations of
su(2) is an important problem, which remains largely unsolved. In principle, a
coproduct can be induced from the ordinary coupling rule for su(2)
\cite{Poly,Cut1} by
the deforming functional technique \cite{Cut2}.
 However such a procedure leads in
general to complicated and untractable coproducts. More direct methods were
recently used to construct a coproduct for some Polychronakos-Ro\v cek
algebras (PRA's) \cite{Grano}
 and for another deformation of su(2) involving a single
function of $J_0$ \cite{Ludu}.\par
%
%
In the present paper, we shall address the problem of endowing some
Delbecq-Quesne algebras (DQA's) with a Hopf algebraic structure. For such
purpose, we shall use a variant and extension of the deforming functional
technique, wherein functionals mapping PRA generators to those of some DQA's
will be determined. By considering the special case where the PRA reduces to
\su, we shall obtain DQA's whose Hopf algebraic structure can be inferred
from that of \su.\par
%
%
This procedure will be carried out in detail for a specific example of DQA, referred
to as \alg, which will be shown to have two sets of ($N+1$)-dimensional unitary
irreducible representations (unirreps), where $N=0$, 1, 2,~$\ldots$. For such an
algebra, the functional to be considered will be two-valued, so the same will be
true for the resulting Hopf algebraic structure. The meaning of this property will
be clarified as far as representation theory is concerned.\par 
%
%
To allow the coupling of any two representations of \alg, we shall be led to
enlarge the double Hopf structure of the algebra into a generalized Hopf structure,
obeying extended coassociativity, counit, and antipode axioms, which we
propose to call a two-colour Hopf algebra. We shall indeed prove that the
generalized coproduct non-cocommutativity is controlled by a generalized
\cR-matrix, satisfying the coloured Yang-Baxter equation 
\cite{Murakami}--\cite{Ohtsuki}, where the `colour' parameters
take two discrete values associated with the two sets of unirreps with the same
dimension that the algebra possesses.\par
%
%
This paper is organized as follows. In Sec.~\ref{sec:functionals} , deforming
functionals mapping PRA's to DQA's are introduced with special emphasis on the
case where the PRA is \su. The algebra \alg is introduced in
Sec.~\ref{sec:example}. In Sec.~\ref{sec:double}, the theory developed in
Sec.~\ref{sec:functionals} is applied to endow \alg\ with a double Hopf algebraic
structure. In Sec.~\ref{sec:generalized}, the latter is enlarged into a generalized
Hopf structure, and the corresponding \cR-matrix is studied in Sec.~\ref{sec:R}.
Sec.~\ref{sec:remarks} contains some concluding remarks.\par
%
%
\section{DEFORMING FUNCTIONALS MAPPING PRA'S TO DQA'S}
\label{sec:functionals}
PRA's are associative algebras over $\C$ generated by three operators
$j_0 = (j_0)^{\dagger}$, $j_+$, and $j_- = (j_+)^{\dagger}$, satisfying the
commutation relations \cite{Poly,Rocek}
\begin{equation}
[j_0, j_+] = j_+, \qquad [j_0, j_-] = - j_-, \qquad [j_+, j_-] = f(j_0),
\label{eq:PRA}
\end{equation} 
where $f(z)$ is a real, parameter-dependent function of $z$, holomorphic in
the neighbourhood of zero, and going to $2z$ for some values of the
parameters. These algebras have a Casimir operator given by
\begin{equation}
c= j_- j_+ + h(j_0) = j_+ j_- + h(j_0) - f(j_0), 
\label{eq:PRA-Casimir}
\end{equation}
in terms of another real function $h(z)$, related to $f(z)$ through the
equation $h(z) - h(z - 1) = f(z)$. In the case where $f(z)$ is a
$\lambda$-degree polynomial, an explicit expression for $h(z)$ has been
found by Delbecq and Quesne~\cite{Quesne}
in terms of Bernoulli polynomials and
Bernoulli numbers.\par
%
%
In the case of polynomial functions $h(z)$, $f(z)$, the PRA can be
identified with the linear space spanned by the monomials
\begin{equation}
  j_-^m j_0^n j_+^p \quad {\rm where} \quad m,\, n,\, p \in \N.
  \label{eq:PRA-basis}
\end{equation}
The above basis is not unique, because we can consider
other bases corresponding to different normal orderings, as the following
ones:
\begin{equation}
  j_+^m j_0^n j_-^p \quad {\rm where} \quad m,\, n,\, p \in \N,
  \label{eq:PRA-basis1}
\end{equation}
or
\begin{equation}
  j_0^m j_+^n j_-^p \quad {\rm where} \quad m,\, n,\, p \in \N.
  \label{eq:PRA-basis2}
\end{equation}
Any normal ordering is actually permitted.\par
%
%
DQA's differ from PRA's by the replacement of~(\ref{eq:PRA}) by \cite{Quesne}
\begin{equation}
[J_0, J_+] = G(J_0) J_+, \qquad [J_0, J_-] = - J_- G(J_0), \qquad
       [J_+, J_-] = F(J_0), 
\label{eq:DQA}
\end{equation}
where the generators $J_0 = (J_0)^{\dagger}$, $J_+$, $J_- = (J_+)^{\dagger}$ are
now denoted by capital letters to distinguish them from those of PRA's,
 and their
commutators involve two real, parameter-dependent functions of $z$,
$F(z)$ and $G(z)$, holomorphic in the neighbourhood of zero, and going to
$2z$ and~1 for some values of the parameters, respectively. These functions
are further restricted by the assumption that the algebras have a Casimir
operator given by 
\begin{equation}
C = J_- J_+ + H(J_0) = J_+ J_- + H(J_0) - F(J_0), 
\label{eq:DQA-Casimir}
\end{equation}
in terms of some real function $H(z)$, holomorphic in the neigh\-bour\-hood 
of zero. The latter restriction 
implies that $F(z)$, $G(z)$, and $H(z)$ satisfy the consistency condition
 $H(z) - H\bigl(z - G(z)\bigr) = F(z)$.\par
%
%
As in the case of the PRA, 
the existence of polynomial functions $G(z)$, $F(z)$, implies
that the DQA can be
identified with the linear space spanned by the monomials
\begin{equation}
  J_-^m J_0^n J_+^p \quad {\rm where} \quad m,\, n,\, p \in \N.
  \label{eq:DQA-basis}
\end{equation}
We must notice here that in general, if the function $G(z)$ is not invertible,
only the normal ordering~(\ref{eq:DQA-basis}) is permitted and 
different orderings, as in the cases~(\ref{eq:PRA-basis1})
and~(\ref{eq:PRA-basis2}), are not allowed.\par
%
%
Let us now try to find a deforming functional that converts the
 PRA generators
into operators 
satisfying the commutation relations of a DQA.
 For such purpose,
we first remark that the first equation in~(\ref{eq:PRA}) can be rewritten as
\begin{equation}
  \bigl(j_0 - 1\bigr) j_+ = j_+ j_0.   \label{eq:prepentire}
\end{equation}
Then, for every entire function $p(z)$, one can prove
\begin{equation}
p\bigl(j_0 - 1\bigr) j_+ = j_+ p(j_0).  \label{eq:entire}
\end{equation}
Let us consider the equation
\begin{equation}
p(z) - p(z-1)= G(p(z)) 
\label{eq:basic}
\end{equation}
for a given function $G(z)$. If this equation has a solution $p(z)$ 
that is an entire function, 
then Eq.~(\ref{eq:entire}) can be written as $[p(j_0),
j_+] = G\left(p(j_0)\right)j_+$. 
Similarly, one finds $[p(j_0), j_-] = - j_-G\left(p(j_0)\right)$.\par
%
%
The above equations suggest that a correspondence between the PRA's and the
DQA's might exist. In terms of the  function $p(z)$, 
the first two relations in~(\ref{eq:PRA})
can indeed be reduced to the first two relations
in~(\ref{eq:DQA}) with the identification
\begin{equation}
J_0 = p(j_0), \qquad J_+ = j_+, \qquad J_- = j_-. 
\label{eq:map}
\end{equation}
Assuming $p$ is invertible, the third defining 
equation~(\ref{eq:PRA}) of a PRA can then be
reduced to the third defining equation~(\ref{eq:DQA})
of a DQA with the identification
\begin{equation}
F \circ p =  f \quad \hbox{\rm or} \quad F = f \circ g, \quad {\rm where} \quad
g=p^{-1},
\label{eq:Ftof}
\end{equation}
and $f\circ g$ means the
composition of the two functions, i.e., $(f\circ g)(z) = f\bigl(g(z)\bigr)$.
The correspondence
between the relevant Casimir operators, 
given in~(\ref{eq:PRA-Casimir}) and~(\ref{eq:DQA-Casimir}), similarly
implies 
\begin{equation}
H \circ p = h \qquad \hbox{\rm or} \qquad H = h \circ g. 
\label{eq:Htoh}      
\end{equation}
\par
%
%
The existence of the map transferring the commutation relations (\ref{eq:PRA}) to
the commutation relations (\ref{eq:DQA}), using the function $p(z)$, does not mean
that the algebra generated by  Eqs.~(\ref{eq:map}) is the same as the DQA. This fact
is evident in the case where the function $g(z)$ is  a polynomial, because the
map~(\ref{eq:map}) tranfers the PRA basis given by  (\ref{eq:PRA-basis}) to the
basis
\begin{equation}
  J_-^m g(J_0)^n J_+^p \quad {\rm where} \quad m,\, n,\, p \in \N,
  \label{eq:map-basis}
\end{equation}
and the latter constitutes a subspace of the linear space~(\ref{eq:DQA-basis})
defining the DQA. Hence, we have shown the following proposition:

\bigskip
{\bf Proposition 1}:
{ \em 
If there exists an invertible entire function $p(z)$ satisfying the equations
$$ p(z) - p(z-1)= G(p(z)) \quad {\rm and} \quad F\circ p =f,$$
then the PRA, defined by the commutation relations
$$[j_0, j_+] = j_+, \quad [j_0, j_-] = - j_-, \quad [j_+, j_-] = f(j_0),$$
is mapped through the generator mapping $P$:
$$J_0 = p(j_0), \qquad J_+ = j_+, \qquad J_- = j_-,$$
into the DQA, defined by the  commutation relations
$$
[J_0, J_+] = G(J_0) J_+, \quad [J_0, J_-] = - J_- G(J_0), \quad
       [J_+, J_-] = F(J_0).$$
} 
\par 
%
%
For a given DQA, the existence of a solution $p(z)$ of Eq.~(\ref{eq:basic}) implies
that the algebraic (and eventually the coalgebraic) structure of the root algebra
PRA can be mapped to the DQA generated by~$P$.\par
%
%
Let consider the special case where the PRA is \su, i.e,
the function $f(j_0)$ in Eq.~(\ref{eq:PRA}) is given by
\begin{equation}
f(j_0) = [2j_0]_q, 
\label{eq:suq2}
\end{equation}
where $[x]_q \equiv (q^x - q^{-x})/(q - q^{-1})$ denotes a $q$-number and $q$
is either a real number or a phase \cite{Biedenharn,Macfarlane}.
Throughout this paper, we shall
assume that $q\in\R^+$. The function $h(j_0)$ of 
Eq.~(\ref{eq:PRA-Casimir}) can be taken as
$h(j_0) = [j_0]_q [j_0 + 1]_q$.\par
%
%
The Hopf algebraic structure of \su~\cite{Drinfeld,Majid} is defined in terms of a
comultiplication map $\Delta: \su \to \su \otimes \su$, a counit map $\epsilon: \su
\to \C$, and an antipode map $S: \su \to \su$, given by
\begin{eqnarray}
   \Delta(j_0) &=& j_0 \otimes 1 + 1 \otimes j_0, \qquad 
       \Delta(j_{\pm}) = j_{\pm} \otimes q^{j_0} + q^{-j_0} \otimes j_{\pm},
       \nonumber \\
   \epsilon(j_0) &=& \epsilon(j_{\pm}) = 0, \nonumber  \\
   S(j_0) &=& - j_0, \qquad S(j_{\pm}) =  - q^{\pm1} j_{\pm}, 
\end{eqnarray}
respectively. These maps satisfy the following relations:
\begin{eqnarray}
   (\Delta \otimes \id) \circ \Delta(a) &=& (\id \otimes \Delta) \circ \Delta(a),
        \nonumber \\
   (\epsilon \otimes \id) \circ \Delta(a) &=& (\id \otimes \epsilon) \circ \Delta(a) =
        a, \nonumber \\
   m \circ (S \otimes \id) \circ \Delta(a) &=& m \circ (\id \otimes S) \circ \Delta(a)
        = \iota \circ \epsilon(a).                   \label{eq:Hopf1}
\end{eqnarray}           
Here $a$ denotes any element of \su, $m$ its multiplication map, $m: \su
\otimes \su \to \su$, and $\iota$ its unit map, $\iota: \C \to \su$, defined by $m(a
\otimes b) = ab$ and $\iota(\lambda) = \lambda 1$ respectively, where 1 is the unit
element of \su. In addition, $\Delta$ and $\epsilon$ are algebra homomorphisms,
\begin{equation}
\Delta(ab) = \Delta(a) \Delta(b), \quad \Delta(1) = 1 \otimes 1, \quad
      \epsilon(ab) = \epsilon(a) \epsilon(b), \quad \epsilon(1) = 1.
\label{eq:Hopf2}
\end{equation} 
\par
%
%
The functions $F(J_0)$ and $H(J_0)$ of Eqs.~(\ref{eq:Ftof}) and~(\ref{eq:Htoh})
become
\begin{eqnarray}
    F(J_0) &=& - {\bigl(\phi(J_0)\bigr)^2 - \bigl(\phi(J_0)\bigr)^{-2}
       \over q - q^{-1}}, \nonumber \\
    H(J_0) &=& {q^{-1} \bigl(\phi(J_0)\bigr)^2 + q \bigl(\phi(J_0)\bigr)^{-2}
       - q - q^{-1}\over (q - q^{-1})^2},    \label{eq:su2-DQA}
\end{eqnarray}
where
$$\phi(z) \equiv q^{-g(z)}. $$
The Hopf algebraic structure of \su\ can then be transferred to that member of
the DQA set whose function $F(J_0)$ is given by~(\ref{eq:su2-DQA}), provided the
function~$g(z)$ can be determined. It should be noticed that the solution of
Eq.~(\ref{eq:basic}) is not a trivial task in general. In the next sections, we shall
proceed to study a special case, which proves tractable.\par
%
%
\section{THE LINEAR $G(J_0)$ CASE} \label{sec:example}
\setcounter{equation}{0}
Let consider the case where the function $G(J_0)$ is linear, i.e.,
\begin{equation}
G(J_0) = 1 + (1-q) J_0, 
\label{eq:linear-G}
\end{equation}
and the $q$ values are restricted to the interval
$0<q<1$ (note that the algebras with $q>1$ are related to those
with $0<q<1$ by an automorphism). The DQA's, corresponding
to~(\ref{eq:linear-G}) and to functions $F(J_0)$ that are polynomials of degree
$\lambda$ in $J_0$, have been studied in Ref.~\cite{Quesne}, where they   
are denoted by ${\cal A}^+_{\alpha_2 \alpha_3 \ldots \alpha_{\lambda-1}
q}(\lambda,1)$, in terms of some extra real parameters $\alpha_2$,~$\alpha_3$, 
$\ldots$,~$\alpha_{\lambda-1}$, entering the
definition of $F(J_0)$. In particular, the representation theory of the algebras
${\cal A}^+_q(2,1)$ (equivalent to Witten's first deformation~\cite{Witten} of
su(2)), and
${\cal A}^+_{p,q}(3,1)$ has been studied in detail. In this respect, the point $z =
(q-1)^{-1}$, where $G(z)$ vanishes, appears as a singular point. The unirreps indeed
separate into two classes according to whether the eigenvalues of $J_0$ are
contained in the interval
$\bigl(-\infty, (q-1)^{-1}\bigr)$, or in the interval $\bigl((q-1)^{-1},
+\infty\bigr)$. In addition, there is a one-dimensional unirrep
corresponding to the eigenvalue $(q-1)^{-1}$ of $J_0$.\par
%
%
With  a linear choice like Eq.~(\ref{eq:linear-G}) for $G(z)$,   
Eq.~(\ref{eq:basic})  can be easily solved. One finds 
a family of  solutions
\begin{equation}
p_\delta (z)=\frac{1- \delta q^{-z} }{q-1},
\label{eq:p-function}
\end{equation} 
where any real, nonvanishing value of the parameter $\delta$ is acceptable.
The functions  $p_\delta$ are entire functions, which are invertible as
\begin{equation}
  p_\delta^{-1}(z)=g_{|\delta|} (z) = \frac{ \ln \left( \left( ( 1 + (1-q)z \right)^2
  /|\delta|^2\right)  }{  \ln \left(1/q^2 \right)}= \frac{ \ln \left( G^2(z) /|\delta|^2
  \right)  }{  \ln \left(1/q^2 \right)}.  \label{eq:inversep} 
\end{equation}
If $z$ is real, the range of $p_{\delta}$ (and consequently the domain of
$p_{\delta}^{-1}$) is the interval
$\bigl(-\infty, (q-1)^{-1}\bigr)$
or $\bigl((q-1)^{-1}, +\infty\bigr)$ according to whether $\delta<0$ or
$\delta>0$. In the case of the linear $G(J_0)$ model, Eq.~(\ref{eq:su2-DQA}) 
gives
\begin{eqnarray}
   F(J_0) &=& - {\bigl(G(J_0)/|\delta|\bigr)^2 - \bigl(G(J_0)/|\delta|\bigr)^{-2} \over 
       q - q^{-1}}, \nonumber \\
   H(J_0) &=& {q^{-1} \bigl(G(J_0)/|\delta|\bigr)^2 + q \bigl(G(J_0)/|\delta|\bigr)^{-2}
      - q - q^{-1}\over (q - q^{-1})^2}.    \label{eq:FG}
\end{eqnarray}
Without loss of generality, we can set $\delta= \pm 1$, because any
$|\delta|\ne 0$ would lead to similar results. Therefore, we shall henceforth
use $p_\pm(z)$ and its inverse $p_\pm^{-1}=g(z)$. 
In this case,
\begin{equation}
p_\pm^{-1}(z)=g(z)= 
\frac{ 
\ln \left( \left( ( 1 + (1-q)z \right)^2 \right) 
}{  \ln \left(1/q^2 \right)}=
\frac{ 
\ln \left( G^2(z) \right) 
}{  \ln \left(1/q^2 \right)}.
\label{eq:g-function}
\end{equation}
\par
%
%
As for any PRA, $j_+$ is the raising generator, while
for the DQA with linear $G(J_0)$, 
$J_+$ (resp.~$J_-$) is the raising generator in the interval 
$\bigl((q-1)^{-1}, +\infty\bigr)$  corresponding to 
$\delta=+1$ (resp.~$\bigl(-\infty, (q-1)^{-1}\bigr)$, corresponding
to $\delta=-1$), one
has to use there the map $p_+(z)$ (resp.~$p_-(z)$).
The function $g(z)$ is well-behaved  everywhere on $\R$, 
except in the neighbourhood of the point $z = (q-1)^{-1}$.\par
%
%
Let denote by \alg\ the DQA generated from \su\ by the
mapping $p_\delta$, with $\delta=\pm 1$:
\begin{equation}
  P_{\delta} : \su \to \alg.       \label{eq:Pdelta}
\end{equation}
\par 
%
%
By a procedure similar to that used in Ref.~\cite{Quesne}, it can be easily shown
that \alg\ has no infinite-dimensional unirrep, but has two ($N+1$)-dimensional
unirreps for any $N = 0$, 1, 2,~$\ldots$. The corresponding spectrum of $J_0$ is
given by
\begin{equation}
 m^{\delta}  = {1 - \delta 
     q^{-(N-2n)/2} \over q-1}, 
\qquad n = 0, 1, \ldots, N, \qquad \delta=\pm 1, 
\label{eq:DQA-J0-eigenvalues}
\end{equation}
with maximum and minimum eigenvalues 
\begin{equation}
J^{\delta} =   {1 -\delta q^{-\delta N/2} \over
    q-1}, \qquad 
-j^{\delta} = {1 -\delta q^{\delta N/2} \over q-1}, 
\label{eq:max-min}
\end{equation}
respectively. The unirrep specified by 
$J^+$ (resp.~$J^-$) is entirely contained
in the interval $\bigl((q-1)^{-1}, +\infty \bigr)$ (resp.~$\bigl(-\infty,
(q-1)^{-1}\bigr)$). For both unirreps, the eigenvalue of the Casimir operator is
given by
\begin{equation}
\langle C\rangle =H(\gamma^\delta),
\label{eq:Casimir}
\end{equation}
where $\gamma^{\delta}= \frac{1- \delta q^{-N/2} }{q-1}$.\par
%
%
In the carrier space $V^{J^{\delta}}$ of the unirrep characterized by~$J^{\delta}$,
whose basis vectors are specified by the values of $J^{\delta}$ and~$m^{\delta}$,
the \alg\ generators are represented by some linear operators
$\Phi^{J^{\delta}}(A)$, $A \in \alg$, defined by
\begin{eqnarray}
   \Phi^{J^{\delta}}\left(J_0\right) \left\vert J^\delta, m^\delta \right> &=&
       m^\delta \left\vert J^\delta, m^\delta \right>, \nonumber\\
   \Phi^{J^{\delta}}\left(J_-\right) \left\vert J^\delta, m^\delta \right> &=&
       \sqrt{ H(\gamma^{\delta}) - H(q m^\delta-1)}
       \left\vert J^\delta, qm^\delta-1 \right>, \nonumber\\
   \Phi^{J^{\delta}}\left(J_+\right) \left\vert J^\delta, m^\delta \right> &=&
       \sqrt{ H(\gamma^{\delta}) - H( m^\delta)}
       \left\vert J^\delta, q^{-1}(m^\delta +1) \right>.   \label{eq:A-modules}
\end{eqnarray}
\par
%
%
\section{DOUBLE HOPF STRUCTURE OF THE ALGEBRA \alg}  \label{sec:double}
\setcounter{equation}{0}
The theory developed in Sec.~\ref{sec:functionals} can be applied provided one
replaces the entire function $p(z)$ by the  functions $p_+(z)$
and $p_-(z)$ 
given by Eq.~(\ref{eq:p-function}), with $\delta=\pm 1$.\par 
%
%
The  functionals $p_\delta$ can be used to transfer the coalgebra
structure and the antipode map from \su to \alg.
The existence of two functionals however leads to similar properties for the latter.
According to whether $p_+$ or $p_-$ is employed, one obtains a comultiplication
map $\Delta_+$ or $\Delta_-$, defined by
\begin{eqnarray}
   \Delta_\delta \left(J_0 \right) &=& \Delta\left( p_\delta (j_0) \right) =
       \case{1}{q-1} \left( 1 \otimes 1 - \delta q^{-j_0} \otimes  q^{-j_0} \right),
       \nonumber \\
   \Delta_\delta \left( J_\pm \right) &=& \Delta \left(j_{\pm}\right) =
       j_\pm \otimes q^{j_0} + q^{-j_0} \otimes j_\pm.   \label{eq:Delta-def}
\end{eqnarray}
The operators $\Delta_\delta \left(J_0 \right)$,
$\Delta_\delta \left( J_\pm \right)$ satisfy the commutation
 relations~(\ref{eq:DQA}) and from the properties of the
comultiplication~(\ref{eq:Hopf1}) and (\ref{eq:Hopf2}),
they satisfy the equation
\begin{equation}
  \Delta_\delta( J_0 ) = p_\delta( \Delta(j_0) ),   \label{eq:Delta-rep}
\end{equation}
hence correspond to a representation of the 
algebra \alg\ characterized by~$\delta$.\par
%
%
The above relations define the comultiplication rules 
as morphisms $\Delta_\delta$ from the algebra \alg\ to the algebra 
$\su \otimes \su$.\par
%
%
By using the generator mapping~$P_{\delta}$ associated with the
functional~$p_{\delta}$, as defined in Proposition 1 and Eq.~(\ref{eq:Pdelta}), we
can go from the algebra $\su \otimes \su$ to the algebra $\alg\otimes\alg$, 
\begin{equation}
  {\rm P}_\delta \otimes {\rm P}_\delta :\;
  \su\otimes \su
  \;\ra \; \alg \otimes \alg.    \label{eq:morphism}
\end{equation}
The above map transfers the Hopf structure of the algebra \su\ to a Hopf-like
structure and we obtain the following proposition:\par
\bigskip
\noindent {\bf Proposition 2}:
{\em The algebra \alg\ is equipped with comultiplication, counit and antipode
maps, given by
\begin{eqnarray}
    \Delta_\delta \left(J_0 \right) &=& \case{1}{q-1} \left(1\otimes 1 - \delta 
        G(J_0)\otimes G(J_0) \right), \nonumber \\
    \Delta_\delta \left( J_\pm \right) &=& \delta\left( J_{\pm} \otimes \left(
        G(J_0) \right)^{-1}  + G(J_0)\otimes J_{\pm} \right), \nonumber\\
    \epsilon_\delta(J_0) &=& \case{1-\delta}{q-1}, \quad 
        \epsilon_\delta (J_\pm) =0, \nonumber\\ 
    S_{\delta}(J_0) &=& - J_0 ( G(J_0) )^{-1}, \quad S_{\delta}(J_+) = - q J_+,
        \quad S_{\delta}(J_-) = - q^{-1} J_-,   \label{eq:Com-def}
\end{eqnarray}
respectively. Both $\Delta_+$, $\epsilon_+$, $S_+$, and
$\Delta_-$, $\epsilon_-$, $S_-$ 
satisfy the Hopf algebra axioms~(\ref{eq:Hopf1}) and~(\ref{eq:Hopf2}),
but the former are only valid for the representations of $\alg\otimes\alg$ 
with eigenvalues of $\Delta_+(J_0)$
in the interval $\bigl((q-1)^{-1}, +\infty \bigr)$,
whereas the latter act in $\bigl(-\infty, (q-1)^{-1}\bigr)$.
 The algebra \alg\ is
therefore endowed with a double Hopf algebraic structure.}
\par\bigskip
\noindent{\it Remark.} Contrary to the comultiplication and counit maps, the
antipode one does not depend explicitly upon $\delta$.\par 
%
%
As a consequence of Eq.~(\ref{eq:Com-def}), from the operators $\varphi^{N/2}(a)$,
$a \in \su$, representing the \su\ generators in the
($N+1$)-dimensional unirrep carrier space $v^{N/2}$, spanned by the vectors
$\left|\frac{N}{2}, \frac{N}{2} - n\right\rangle$, $n=0$,~1, $\ldots$,~$N$,
\cite{Biedenharn,Macfarlane}
\begin{eqnarray}
   \varphi^{N/2}(j_0) \left|\case{N}{2}, \case{N}{2} - n\right\rangle &=&
        \left(\case{N}{2} - n\right) \left|\case{N}{2}, \case{N}{2} - n\right\rangle,
        \nonumber \\
   \varphi^{N/2}(j_-) \left|\case{N}{2}, \case{N}{2} - n\right\rangle &=&
        \sqrt{[n+1]_q [N-n]_q} \, \left|\case{N}{2}, \case{N}{2} - n - 1\right\rangle,
        \nonumber \\
   \varphi^{N/2}(j_+) \left|\case{N}{2}, \case{N}{2} - n\right\rangle &=& \sqrt{[n]_q
        [N-n+1]_q} \, \left|\case{N}{2}, \case{N}{2} - n + 1\right\rangle,  \label{eq:j}
\end{eqnarray}
one obtains the operators $\Phi^{J^{\delta}}(A)$, $A \in \alg$, representing the \alg
\ generators in $V^{j^{\delta}}$ as follows:
\begin{eqnarray}
   \Phi^{J^{\delta}}\left(J_0\right) \left|J^{\delta}, m^{\delta}\right\rangle &=&
          m^{\delta} \left|J^{\delta}, m^{\delta}\right\rangle, \nonumber \\
   \Phi^{J^{\delta}}\left(J_-\right) \left|J^{\delta}, m^{\delta}\right\rangle &=&
          \sqrt{[n+1]_q [N-n]_q}\, \left|J^{\delta}, q m^{\delta} - 1\right\rangle,
          \nonumber \\
   \Phi^{J^{\delta}}\left(J_+\right) \left|J^{\delta}, m^{\delta}\right\rangle &=&
          \sqrt{[n]_q [N-n+1]_q}\, \left|J^{\delta}, q^{-1}
          \left(m^{\delta}+1\right)\right\rangle.  \label{eq:J}
\end{eqnarray}
On the right-hand side of~(\ref{eq:J}), $N$ and $n$ have to be replaced by their
expression in terms of $J^{\delta}$ and $m^{\delta}$, 
obtained by inverting~(\ref{eq:DQA-J0-eigenvalues})
and~(\ref{eq:max-min}). The results can be written as 
in Eq.~(\ref{eq:A-modules}), which was constructed
by a direct procedure. They also confirm that \alg\ has no infinite-dimensional
unirrep.\par
%
%
Considering now an $(N_1+1)$-dimensional unirrep of \alg, characterized by
$J^\delta_1$, in a carrier space $V^{J_1^{\delta}}$, and another
$(N_2+1)$-dimensional unirrep of the same, specified by $J^\delta_2$, in a carrier
space $V^{J_2^{\delta}}$, one can couple them by using the coproduct
$\Delta_\delta$  to obtain a reducible representation of \alg\ in $V^{J_1^{\delta}}
\otimes V^{J_2^{\delta}}$. The corresponding operators $\Phi^{J_1^{\delta}
J_2^{\delta}}(A)$, $A \in \alg$, are given by
\pagebreak
\begin{eqnarray}
  \lefteqn{\Phi^{J_1^{\delta} J_2^{\delta}}\left(J_0\right) \left(\left|J_1^{\delta},
       m_1^{\delta}\right\rangle \otimes \left|J_2^{\delta}, m_2^{\delta}
       \right\rangle\right) } \nonumber \\ 
  & & = \case{1}{q-1} \left(1 - \delta G\left(m_1^{\delta}
       \right) G\left(m_2^{\delta}\right)\right) \left|J_1^{\delta},
       m_1^{\delta}\right\rangle \otimes \left|J_2^{\delta}, m_2^{\delta}
       \right\rangle, \nonumber\\
  \lefteqn{\Phi^{J_1^{\delta} J_2^{\delta}}\left(J_-\right) \left(\left|J_1^{\delta},
       m_1^{\delta}\right\rangle \otimes \left|J_2^{\delta}, m_2^{\delta}
       \right\rangle\right)} \nonumber \\
  & & = \sqrt{[n_1+1]_q [N_1-n_1]_q}\, q^{(N_2-2n_2)/2} \left|J_1^{\delta},
       q m_1^{\delta}-1\right\rangle \otimes \left|J_2^{\delta},
       m_2^{\delta}\right\rangle \nonumber \\
  & & \phantom{=} \mbox{} + q^{-(N_1-2n_1)/2} \sqrt{[n_2+1]_q [N_2-n_2]_q}\, 
       \left|J_1^{\delta}, m_1^{\delta}\right\rangle \otimes \left|J_2^{\delta},
       q m_2^{\delta}-1\right\rangle, \nonumber \\ 
  \lefteqn{\Phi^{J_1^{\delta} J_2^{\delta}}\left(J_+\right) \left(\left|J_1^{\delta},
       m_1^{\delta}\right\rangle \otimes \left|J_2^{\delta}, m_2^{\delta}
       \right\rangle\right)} \nonumber \\
  & & = \sqrt{[n_1]_q [N_1-n_1+1]_q}\, q^{(N_2-2n_2)/2} \left|J_1^{\delta},
       q^{-1} \left(m_1^{\delta}+1\right)\right\rangle \otimes \left|J_2^{\delta},
       m_2^{\delta}\right\rangle \nonumber \\
  & & \phantom{=} \mbox{} + q^{-(N_1-2n_1)/2} \sqrt{[n_2]_q [N_2-n_2+1]_q}
       \nonumber \\ 
  & & \phantom{=} \mbox{} \times\left|J_1^{\delta}, m_1^{\delta}\right\rangle
       \otimes \left|J_2^{\delta}, q^{-1}\left(m_2^{\delta}+1\right)\right\rangle,
  \label{eq:reducible} 
\end{eqnarray}
where $m^{\delta}_1$, $m^{\delta}_2$ are defined in terms of $n_1$, $n_2$,
respectively,  in the same way as $m^{\delta}$ in terms of $n$ (see
Eq.~(\ref{eq:DQA-J0-eigenvalues})).\par
%
%
Such a reducible representation can be decomposed into a direct sum of
$(N+1)$-dimensional unirreps, characterized by $J^\delta$, and whose basis states
$\left| J_1^{\delta} J_2^{\delta} J^{\delta} m^{\delta} \right\rangle$ can be
written as
\begin{equation}
  \left| J_1^{\delta} J_2^{\delta} J^{\delta} m^{\delta} \right\rangle = 
  \sum_{m_1^{\delta}, m_2^{\delta}} \left\langle J^{\delta}_1\, m^{\delta}_1,  
  J^{\delta}_2\, m^{\delta}_2 \right. \left| J^{\delta}\, m^{\delta}
  \right\rangle_{DQ} \left|J_1^{\delta}, m_1^{\delta}\right\rangle \otimes
  \left|J_2^{\delta}, m_2^{\delta} \right\rangle,   \label{eq:coupled}  
\end{equation}
in terms of some Wigner coefficients $\langle J^\delta_1\, m^\delta_1, 
J^\delta_2\, m^\delta_2 | J^\delta\, m^\delta \rangle_{DQ}$. From
the relation between the representations of \alg\ and those of \su, it
follows that
\begin{equation}
   \left\langle J^{\delta}_1\, m^{\delta}_1,  J^{\delta}_2\, m^{\delta}_2 \right.
   \left| J^{\delta}\, m^{\delta} \right\rangle_{DQ} = 
   \left\langle \case{N_1}{2}\; \case{N_1}{2} - n_1, 
   \case{N_2}{2}\; \case{N_2}{2} - n_2 \right. \left| \case{N}{2}\; \case{N}{2} - n
   \right\rangle_q,     \label{eq:Wigner}
\end{equation}
where the quantity on the right-hand side is an \su Wigner
coefficient~\cite{Nomura}. It follows from the  first formula in
Eq.~(\ref{eq:Com-def}) that the Wigner coefficient~(\ref{eq:Wigner}) can be
different from zero only if
\begin{equation}
  m^\delta = \case{1}{q-1} \left( 1- \delta  G(m_1^\delta) G(m_2^\delta) \right).
  \label{eq:totalm}
\end{equation}
\par 
%
%
It is worth noting that up to now no comultiplication rule is available for
coupling two unirreps $J^+_1$ and $J^-_2$, or $J^-_1$ and $J^+_2$. The purpose of
the next section will be to show that the comultiplication rule definition can be
extended so as to allow the coupling of any two \alg\ unirreps. This extension will
lead us to enlarge the double Hopf structure of the algebra into a generalized one.
\par
%
%
\section{GENERALIZED HOPF STRUCTURE OF THE ALGEBRA \alg}   
\label{sec:generalized}
\setcounter{equation}{0}
To connect the two types of unirreps specified by $\delta = +1$ and~$\delta = -1$
respectively, let us introduce some linear operators $T^{J^{\delta}}: V^{J^{\delta}}
\to V^{J^{-\delta}}$, where $J^{\delta}$ may be any unirrep label, as given in
Eq.~(\ref{eq:max-min}). They are defined by their action on the $V^{J^{\delta}}$
basis vectors $\left\vert J^\delta, m^\delta \right>$ as follows:
\begin{equation}
  T^{J^{\delta}} \left\vert J^\delta, m^\delta \right> = \left\vert J^{-\delta},
  m^{-\delta}\right>.     \label{eq:transmutation-def}
\end{equation}
Such operators will be referred to as transmutation operators as they change the
basis states of an ($N+1$)-dimensional unirrep into those of its companion with
the same dimension. They obviously satisfy the relation
\begin{equation}
  T^{J^{-\delta}} T^{J^{\delta}} = I^{J^{\delta}},    \label{eq:transmutation-inv}
\end{equation}
where $I^{J^{\delta}}$ denotes the unit operator in $V^{J^{\delta}}$.\par
%
%
By applying $T^{J^{\delta}}$ on both sides of Eq.~(\ref{eq:J}) and using
Eqs.~(\ref{eq:transmutation-def}), (\ref{eq:transmutation-inv}), and
(\ref{eq:DQA-J0-eigenvalues}), it can be easily proved that for any \alg\
generator~$A$,
\begin{equation}
  T^{J^{\delta}} \Phi^{J^{\delta}}(A)\, T^{J^{-\delta}} = \Phi^{J^{-\delta}}(\sigma(A)),
  \label{eq:transmutation-irrep}
\end{equation}
where $\sigma: \alg \to \alg$, defined by
\begin{equation}
  \sigma(J_0) = \case{2}{q-1} - J_0, \qquad \sigma(J_{\pm}) = J_{\pm},
  \label{eq:sigma}
\end{equation}
is an involutive automorphism of the algebra \alg. This clearly shows that at the
algebra level, the operator~$\sigma$ is responsible for the transmutation.\par
%
%
Let us define $\sigma_{\delta}: \alg \to \alg$ by
\begin{equation}
  \sigma_{\delta} = \left\{\begin{array}{ll}
                                             \id & \mbox{if $\delta=+1$}\\
                                             \sigma & \mbox{if $\delta=-1$}  
                                          \end{array}
                               \right..        \label{eq:sigmadelta}
\end{equation}
The basic mapping~$P_{\delta}$ of Eq.~(\ref{eq:Pdelta}) transforms
under~(\ref{eq:sigmadelta}) as follows:
\begin{equation}
  \sigma_{\zeta\eta} \circ P_{\eta} = P_{\zeta},      \label{eq:sigmaP}
\end{equation}
where $\zeta$, $\eta=\pm 1$. This equation is equivalent to the commuting diagram
\begin{equation}
  \begin{array}{ccc}
  \su & \stackrel{P_{\eta}}{\ra} & \alg \\[0.2cm] 
  {\ss P_{\zeta}}{\da} &  & {\da} {\ss\sigma_{\zeta\eta}} \\[0.2cm]
  \alg & \stackrel{\ids}{\lra} & \alg 
\end{array}     \label{eq:diagsigmaP}
\end{equation}
\par
%
%
We can now extend the comultiplication and antipode maps, $\Delta_{\delta}$ and
$S_{\delta}$ ($\delta = \pm1$), of Eq.~(\ref{eq:Com-def}) by setting
\begin{equation}
  \Delta^{\zeta,\eta}_{\delta}(A) = \left(\sigma_{\zeta \delta} \otimes
  \sigma_{\eta \delta}\right) \circ \Delta_{\delta}(A), \qquad
  S^{\zeta}_{\delta}(A) = \sigma_{\zeta \delta} \circ S_{\delta}(A),   
  \label{eq:extcom-def}
\end{equation}
where $\zeta$, $\eta$, $\delta = \pm1$, while leaving unchanged the counit map
$\epsilon_{\delta}$, defined in the same equation. We note that in particular,
\begin{equation}
  \Delta^{\delta,\delta}_{\delta} = \Delta_{\delta}, \qquad 
  S^{\delta}_{\delta} = S_{\delta}.         \label{eq:oldcom}
\end{equation}
By using Eqs.~(\ref{eq:Com-def}), (\ref{eq:sigma}), and~(\ref{eq:sigmadelta}), we
obtain
\begin{eqnarray}
    \Delta^{\zeta,\eta}_\delta \left(J_0 \right) &=& \case{1}{q-1} \left(1\otimes 1
        - \delta \zeta \eta G(J_0)\otimes G(J_0) \right), \nonumber \\
    \Delta^{\zeta,\eta}_\delta \left( J_\pm \right) &=& \eta J_{\pm} \otimes
        \left(G(J_0) \right)^{-1}  + \zeta G(J_0)\otimes J_{\pm},
        \nonumber\\
    S^{\zeta}_{\delta}(J_0) &=& \case{1}{q-1} \left(1 - \zeta \delta
        \left(G(J_0)\right)^{-1}\right), \nonumber \\ 
    S^{\zeta}_{\delta}(J_{\pm}) & = & - q^{\pm1} J_{\pm}.      \label{eq:extcom}
\end{eqnarray}
\par
%
%
Alternatively, $\Delta^{\zeta,\eta}_{\delta}$, $\epsilon_{\delta}$, and
$S^{\zeta}_{\delta}$ can be defined directly in terms of the comultiplication, counit,
and antipode maps $\Delta$, $\epsilon$, $S$ of \su, as well as the
map~$P_{\delta}$, by the commuting diagrams
\begin{equation}
  \begin{array}{ccc}
  \su & \stackrel{\Delta}{\ra} & \su\otimes\su \\[0.2cm] 
  {\ss P_{\delta}}{\da} &  & {\da} {\ss P_{\zeta}\otimes P_{\eta}} \\[0.2cm]
  \alg & \stackrel{\Delta^{\zeta,\eta}_{\delta}}{\ra} & \alg\otimes\alg 
\end{array}     
\qquad
  \begin{array}{ccc}
  \su & \stackrel{\epsilon}{\ra} & \C \\[0.2cm] 
  {\ss P_{\delta}}{\da} &  & {\uda} {\ids} \\[0.2cm]
  \alg & \stackrel{\epsilon_{\delta}}{\ra} & \C 
\end{array}
\qquad 
  \begin{array}{ccc}
  \su & \stackrel{S}{\ra} & \su \\[0.2cm] 
  {\ss P_{\delta}}{\da} &  & {\da} {\ss P_{\zeta}} \\[0.2cm]
  \alg & \stackrel{S^{\zeta}_{\delta}}{\ra} & \alg 
\end{array}
\label{eq:diagDelta}
\end{equation}
\par
%
%
As shown in the Appendix, $\Delta^{\zeta,\eta}_{\delta}$, $\epsilon_{\delta}$, and
$S^{\zeta}_{\delta}$ transform under $\sigma_{\delta}$ as follows:
\begin{eqnarray}
  \left(\sigma_{\mu \zeta} \otimes \sigma_{\nu \eta}\right) \circ
         \Delta^{\zeta,\eta}_{\delta} & = & \Delta^{\mu,\nu}_{\rho} \circ
         \sigma_{\rho\delta}, \label{eq:extcom-transf1} \\
  \epsilon_{\delta} \circ \sigma_{\delta \zeta} & = & \epsilon_{\zeta},
         \label{eq:extcom-transf2} \\
  \sigma_{\zeta \eta} \circ S^{\eta}_{\delta} & = & S^{\zeta}_{\mu} \circ
         \sigma_{\mu \delta}.   
         \label{eq:extcom-transf3}
\end{eqnarray}
\par
%
%
By using Eqs.~(\ref{eq:extcom-transf1})--(\ref{eq:extcom-transf3}) and the
Hopf algebra axioms  (\ref{eq:Hopf1}),~(\ref{eq:Hopf2}), satisfied by
$\Delta_{\delta}$,
$\epsilon_{\delta}$, and~$S_{\delta}$, or alternatively the diagrammatic method of
the Appendix, we obtain
\par
\bigskip
\noindent {\bf Proposition 3}:
{\em The algebra \alg\ is endowed with a generalized Hopf algebraic structure,
whose comultiplication, counit, and antipode maps, $\Delta^{\zeta,\eta}_{\delta}$,
$\epsilon_{\delta}$, $S^{\zeta}_{\delta}$, defined in Eqs.~(\ref{eq:extcom}) and
(\ref{eq:Com-def}), satisfy the following generalized coassociativity, counit, and
antipode axioms: 
\begin{eqnarray}
  \left(\Delta^{\zeta,\eta}_{\mu} \otimes \id\right) \circ \Delta^{\mu,\nu}_{\delta}
        (A) & = & \left(\id \otimes \Delta^{\eta,\nu}_{\rho}\right) \circ
        \Delta^{\zeta,\rho}_{\delta}(A), \label{eq:genHopf1} \\
  \left(\epsilon_{\zeta} \otimes \sigma_{\eta \delta}\right) \circ
        \Delta^{\zeta,\eta}_{\delta}(A) & = & \left(\sigma_{\zeta \delta} \otimes
        \epsilon_{\eta}\right) \circ \Delta^{\zeta,\eta}_{\delta}(A) = A, 
        \label{eq:genHopf2}\\
  m \circ \left(S^{\mu}_{\zeta} \otimes \sigma_{\mu \eta}\right) \circ
        \Delta^{\zeta,\eta}_{\delta}(A) & = & m \circ \left(\sigma_{\mu \zeta}
        \otimes S^{\mu}_{\eta}\right) \circ \Delta^{\zeta,\eta}_{\delta}(A)
        \nonumber \\ 
        & = & \iota \circ \epsilon_{\delta}(A),      \label{eq:genHopf3}
\end{eqnarray}
where $A$ denotes any element of \alg, $m$ and $\iota$ are the multiplication and
unit maps of \alg, $\delta$, $\zeta$, $\eta$, $\mu$, $\nu$, $\rho$ take any values in
the set $\{-1,+1\}$, and no summation over repeated indices is implied. Moreover,
$\Delta^{\zeta,\eta}_{\delta}$ and
$\epsilon_{\delta}$ are algebra homomorphisms, while $S^{\zeta}_{\delta}$ is both
an algebra and a coalgebra antihomomorphism.}
\par\bigskip
\noindent{\it Remark.} In principle, the multiplication and unit maps of 
\alg\ should be distinguished from those of \su, but as no confusion can arise, for
simplicity's sake we use the same notation.\par
%
%
By using the generalized coproduct $\Delta^{\zeta,\eta}_{\delta}$, any ($N_1+1$)-
and ($N_2+1$)-dimensional unirreps of \alg, specified by $J_1^{\zeta}$ and
$J_2^{\eta}$ respectively, can now be coupled to provide two reducible
representations in $V^{J_1^{\zeta}} \otimes V^{J_2^{\eta}}$, which are
characterized by $\delta = +1$ and $\delta = -1$, respectively. From
Eq.~(\ref{eq:extcom}), we obtain for the corresponding operators $\Phi^{J_1^{\zeta}
J_2^{\eta}}_{\delta}(A)$, $A \in \alg$,
\begin{eqnarray}
  \lefteqn{\Phi^{J_1^{\zeta} J_2^{\eta}}_{\delta}\left(J_0\right)
       \left(\left|J_1^{\zeta}, m_1^{\zeta}\right\rangle \otimes \left|J_2^{\eta},
       m_2^{\eta}\right\rangle\right) } \nonumber \\ 
  & & = \case{1}{q-1} \left(1 - \delta q^{-(N_1+N_2-2n_1-2n_2)/2}\right)
       \left|J_1^{\zeta}, m_1^{\zeta}\right\rangle \otimes \left|J_2^{\eta},
       m_2^{\eta}\right\rangle, \nonumber\\
  \lefteqn{\Phi^{J_1^{\zeta} J_2^{\eta}}_{\delta}\left(J_-\right)
       \left(\left|J_1^{\zeta}, m_1^{\zeta}\right\rangle \otimes \left|J_2^{\eta},
       m_2^{\eta} \right\rangle\right)} \nonumber \\
  & & = \sqrt{[n_1+1]_q [N_1-n_1]_q}\, q^{(N_2-2n_2)/2} \left|J_1^{\zeta},
       q m_1^{\zeta}-1\right\rangle \otimes \left|J_2^{\eta},
       m_2^{\eta}\right\rangle \nonumber \\
  & & \phantom{=} \mbox{} + q^{-(N_1-2n_1)/2} \sqrt{[n_2+1]_q [N_2-n_2]_q}\, 
       \left|J_1^{\zeta}, m_1^{\zeta}\right\rangle \otimes \left|J_2^{\eta},
       q m_2^{\eta}-1\right\rangle, \nonumber \\ 
  \lefteqn{\Phi^{J_1^{\zeta} J_2^{\eta}}_{\delta}\left(J_+\right)
       \left(\left|J_1^{\zeta}, m_1^{\zeta}\right\rangle \otimes \left|J_2^{\eta},
       m_2^{\eta}\right\rangle\right)} \nonumber \\
  & & = \sqrt{[n_1]_q [N_1-n_1+1]_q}\, q^{(N_2-2n_2)/2} \left|J_1^{\zeta},
       q^{-1} \left(m_1^{\zeta}+1\right)\right\rangle \otimes \left|J_2^{\eta},
       m_2^{\eta}\right\rangle \nonumber \\
  & & \phantom{=} \mbox{} + q^{-(N_1-2n_1)/2} \sqrt{[n_2]_q [N_2-n_2+1]_q}
       \nonumber \\ 
  & & \phantom{=} \mbox{} \times\left|J_1^{\zeta}, m_1^{\zeta}\right\rangle
       \otimes \left|J_2^{\eta}, q^{-1}\left(m_2^{\eta}+1\right)\right\rangle.
  \label{eq:reduciblebis} 
\end{eqnarray}
%
%
By comparing Eq.~(\ref{eq:reduciblebis}) with Eq.~(\ref{eq:reducible}), we note that
the $\delta$-type reducible representation in $V^{J_1^{\zeta}} \otimes
V^{J_2^{\eta}}$ coincides with the reducible representation in
$V^{J_1^{\delta}} \otimes V^{J_2^{\delta}}$, previously considered. Therefore the
same transformation decomposes both representations into a direct sum of
($N+1$)-dimensional unirreps characterized by~$J^{\delta}$. Hence, the states
\begin{equation}
  \left| J_1^{\zeta} J_2^{\eta} J^{\delta} m^{\delta} \right\rangle = 
  \sum_{m_1^{\zeta}, m_2^{\eta}} \left\langle J^{\zeta}_1\, m^{\zeta}_1,  
  J^{\eta}_2\, m^{\eta}_2 \right. \left| J^{\delta}\, m^{\delta}
  \right\rangle_{DQ} \left|J_1^{\zeta}, m_1^{\zeta}\right\rangle \otimes
  \left|J_2^{\eta}, m_2^{\eta} \right\rangle,   \label{eq:coupledbis}  
\end{equation}
with
\begin{equation}
   \left\langle J^{\zeta}_1\, m^{\zeta}_1,  J^{\eta}_2\, m^{\eta}_2 \right.
   \left| J^{\delta}\, m^{\delta} \right\rangle_{DQ} = 
   \left\langle \case{N_1}{2}\; \case{N_1}{2} - n_1, 
   \case{N_2}{2}\; \case{N_2}{2} - n_2 \right. \left| \case{N}{2}\; \case{N}{2} - n
   \right\rangle_q,     \label{eq:Wignerbis}
\end{equation}
span the carrier space of the unirrep~$J^{\delta}$ in $V^{J_1^{\zeta}} \otimes
V^{J_2^{\eta}}$.\par
%
%
It should be stressed that the space $V^{J_1^{\zeta}} \otimes V^{J_2^{\eta}}$
{\it does not contain two (N+1)-dimensional unirreps}, characterized by~$J^+$
and~$J^-$ respectively, for $N$ in the range $|N_1 - N_2|$, $|N_1 - N_2| + 2$,
$\ldots$,~$N_1 + N_2$, {\it but a single one}, which may be considered as 
that specified by~$J^+$ or that specified by~$J^-$, according to whether the
coproduct $\Delta^{\zeta,\eta}_+$ or $\Delta^{\zeta,\eta}_-$ is used. In other
words, a given linear combination of states $\left|J_1^{\zeta}, m_1^{\zeta}
\right\rangle \otimes \left|J_2^{\eta}, m_2^{\eta} \right\rangle$, as that
contained in Eq.~(\ref{eq:coupledbis}), may be regarded as a basis state of a
unirrep $J^+$ {\it or\/} $J^-$, where $J^+$ and $J^-$ are determined from~$N$ by
using Eq.~(\ref{eq:max-min}) (and in the same way, the corresponding $m^+$ or
$m^-$ is determined from~$N$ and~$n$ by using
Eq.~(\ref{eq:DQA-J0-eigenvalues})).\par 
%
%
For instance, the state $\left|J_1^+, J_1^+\right\rangle \otimes \left|J_2^+, J_2^+
\right\rangle$, where $J_1^+ = J_2^+ =
\left(\sqrt{q}\right.$ $\left.\left(\sqrt{q}+1\right)\right)^{-1}$, corresponding to
$N_1 = N_2 =1$, $n_1 = n_2 = 0$, may be considered as the highest-weight state of
the three-dimensional unirrep characterized by $J^+ = q^{-1}$, or the lowest-weight
state of the three-dimensional unirrep specified by $J^- = (q+1)/(q-1)$, both of
these states corresponding to $N=2$ and $n=0$ in
Eqs.~(\ref{eq:DQA-J0-eigenvalues}) and~(\ref{eq:max-min}),
\begin{eqnarray}
  \left|J_1^+, J_1^+\right\rangle \otimes \left|J_2^+, J_2^+ \right\rangle 
  & = & \left|J_1^+ J_2^+ J^+=q^{-1}, m^+=q^{-1}\right\rangle\nonumber \\
  & = & \left|J_1^+ J_2^+ J^-=\case{q+1}{q-1}, m^-=\case{q+1}{q(q-1)}\right\rangle.
  \label{eq:example-state}
\end{eqnarray}
By direct use of Eq.~(\ref{eq:extcom}), we indeed obtain for instance
\begin{eqnarray}
  \Delta^{+,+}_+\left(J_0\right) \left(\left|J_1^+, J_1^+\right\rangle \otimes
         \left|J_2^+, J_2^+\right\rangle\right) & = & q^{-1}
         \left(\left|J_1^+, J_1^+\right\rangle \otimes
         \left|J_2^+, J_2^+\right\rangle\right), \nonumber \\
  \Delta^{+,+}_-\left(J_0\right) \left(\left|J_1^+, J_1^+\right\rangle \otimes
         \left|J_2^+, J_2^+\right\rangle\right) & = & \frac{q+1}{q(q-1)}
         \left(\left|J_1^+, J_1^+\right\rangle \otimes
         \left|J_2^+, J_2^+\right\rangle\right), \nonumber \\ 
  \Delta^{+,+}_{\pm}\left(J_+\right) \left(\left|J_1^+, J_1^+\right\rangle \otimes
         \left|J_2^+, J_2^+\right\rangle\right) & = & 0.    \label{eq:example-com}
\end{eqnarray}
\par
%
%
In the next section, we shall examine how the universal \cR-matrix definition
valid for \su, hence for the double Hopf structure of \alg, can be extended to the
generalized Hopf structure of the latter.\par
%
%
\section{GENERALIZED $\cal R$-MATRIX OF THE ALGE\-BRA \alg}
\label{sec:R}
\setcounter{equation}{0}
It is well known~\cite{Majid} that \su\ is a quasitriangular Hopf algebra, which
means that there exists an invertible element $\cR \in \su \otimes \su$
(completed tensor product), called the \su\ universal \cR-matrix,
\begin{equation}
  \cR = q^{2j_0 \otimes j_0} \sum_{n=0}^{\infty} \frac{(1-q^{-2})^n}{[n]_q!}\,
  q^{n(n-1)/2} \left(q^{j_0} j_+ \otimes q^{-j_0} j_-\right)^n,     \label{eq:R-su}
\end{equation}
such that its comultiplication~$\Delta$ and its opposite comultiplication
$\Delta^{\mbox{\rm\scriptsize op}} \equiv \tau \circ \Delta$ (where $\tau$ is the
twist map, defined by $\tau(a \otimes b) = b \otimes a$) only differ by a
conjugation by
\cR,
\begin{equation}
  \Delta^{\mbox{\rm\scriptsize op}}(a) = \cR \Delta(a) \cR^{-1}, \qquad a \in \su,
  \label{eq:Rcond1-su}  
\end{equation}
and in addition
\begin{equation}
  (\Delta \otimes \id)(\cR) = \cR_{13} \cR_{23}, \qquad (\id \otimes \Delta)(\cR)
  = \cR_{13} \cR_{12}.     \label{eq:Rcond2-su} 
\end{equation}
\par
%
%
By applying the map $P_{\delta} \otimes P_{\delta}$, considered in
Eq.~(\ref{eq:morphism}), to the \su\ \cR-matrix, given in Eq.~(\ref{eq:R-su}), we
obtain an invertible element $\cR^{\delta}$ of $\alg \otimes \alg$ (completed
tensor product),
\begin{equation}
  \cR^{\delta} = q^{2\log_q(\delta G(J_0)) \otimes  \log_q(\delta G(J_0))}
  \sum_{n=0}^{\infty} \frac{(1-q^{-2})^n}{[n]_q!}\, q^{n(n-1)/2} \left((G(J_0))^{-1}
  J_+ \otimes G(J_0) J_-\right)^n.         \label{eq:Rdelta}
\end{equation}
It can be easily checked that it satisfies with respect to the comultiplication
$\Delta_{\delta}$, defined in Eq.~(\ref{eq:Com-def}), and the opposite
comultiplication $\Delta^{\mbox{\rm\scriptsize op}}_{\delta} \equiv \tau \circ
\Delta_{\delta}$, relations similar to Eqs.~(\ref{eq:Rcond1-su})
and~(\ref{eq:Rcond2-su}), namely
\begin{eqnarray}
  \Delta^{\mbox{\rm\scriptsize op}}_{\delta}(A) & = & \cR^{\delta}
          \Delta_{\delta}(A) \left(\cR^{\delta}\right)^{-1}, \qquad A \in \alg,
          \label{eq:Rdelta-cond1} \\
  (\Delta_{\delta} \otimes \id)\left(\cR^{\delta}\right) & = & \cR^{\delta}_{13}
          \cR^{\delta}_{23}, \qquad (\id \otimes \Delta_{\delta})\left(\cR^{\delta}
          \right) = \cR^{\delta}_{13} \cR^{\delta}_{12}.    \label{eq:Rdelta-cond2}
\end{eqnarray}
Hence, \alg\ has a double quasitriangular Hopf structure with a double universal
\cR-matrix $\cR^{\delta}$, $\delta = \pm1$, where $\cR^+$ (resp.~$\cR^-$)
corresponds to the coalgebra structure and antipode ($\Delta_+$, $\epsilon_+$,
$S_+$) (resp.\  ($\Delta_-$, $\epsilon_-$, $S_-$)), and therefore acts in the
interval $\bigl((q-1)^{-1},+\infty\bigr)$ (resp.\ $\bigl(-\infty,
(q-1)^{-1}\bigr)$).\par
%
%
It is worth noting that as direct consequences of Eqs.~(\ref{eq:Rdelta-cond1}),
(\ref{eq:Rdelta-cond2}), and of the Hopf algebra axioms
(\ref{eq:Hopf1}),~(\ref{eq:Hopf2}), satisfied by $\Delta_{\delta}$,
$\epsilon_{\delta}$,~$S_{\delta}$, the double \cR-matrix $\cR^{\delta}$, $\delta =
\pm1$, satisfies the relations
\begin{eqnarray}
  \cR^{\delta}_{12} \cR^{\delta}_{13} \cR^{\delta}_{23} & = & \cR^{\delta}_{23}
         \cR^{\delta}_{13} \cR^{\delta}_{12},      \label{eq:Rdelta-prop1} \\
  (\epsilon_{\delta} \otimes \id) \left(\cR^{\delta}\right) & = & (\id \otimes
         \epsilon_{\delta}) \left(\cR^{\delta}\right) = 1,     \label{eq:Rdelta-prop2} \\
  (S_{\delta} \otimes \id) \left(\cR^{\delta}\right) & = & \left(\id \otimes
         S_{\delta}^{-1}\right) \left(\cR^{\delta}\right) =
         \left(\cR^{\delta}\right)^{-1},           \label{eq:Rdelta-prop3}
\end{eqnarray}
which are similar to well-known properties of the \su\ \cR-matrix~\cite{Majid}.
In particular, Eq.~(\ref{eq:Rdelta-prop1}) shows that both $\cR^+$ and $\cR^-$ are
solutions of the (ordinary) Yang-Baxter equation (YBE).\par
%
%
Turning now to the generalized Hopf structure introduced in the previous section,
let us consider in $\alg \otimes \alg$ (completed tensor product) the four elements
\begin{equation}
  \cR^{\zeta,\eta} = (\sigma_{\zeta\delta} \otimes \sigma_{\eta\delta}) 
  \left(\cR^{\delta}\right), \qquad \zeta, \eta = \pm1,      \label{eq:genR}
\end{equation}
where, on the right-hand side, $\delta$ takes any value in the set $\{-1,+1\}$
and no summation over it is implied. From Eqs.~(\ref{eq:sigma}),
(\ref{eq:sigmadelta}), and~(\ref{eq:Rdelta}), it follows that the explicit form of
$\cR^{\zeta,\eta}$ is given by
\begin{equation}
  \cR^{\zeta,\eta} = q^{2\log_q(\zeta G(J_0)) \otimes  \log_q(\eta G(J_0))}
  \sum_{n=0}^{\infty} \frac{(1-q^{-2})^n}{[n]_q!}\, q^{n(n-1)/2}
  \left((\zeta G(J_0))^{-1} J_+ \otimes \eta G(J_0) J_-\right)^n.        
  \label{eq:genR-exp}
\end{equation}
We also note that $\cR^{\delta,\delta} = \cR^{\delta}$.\par
%
%
By using now Eqs.~(\ref{eq:extcom-def}), 
(\ref{eq:extcom-transf1})--(\ref{eq:extcom-transf3}), (\ref{eq:Rdelta-cond1}),
(\ref{eq:Rdelta-cond2}), and~(\ref{eq:genR}), we easily obtain
\par
\bigskip
\noindent {\bf Proposition 4}:
{\em The generalized Hopf algebra, defined in Proposition~3 of the previous
section, has a generalized universal \cR-matrix made of four invertible pieces
$\cR^{\zeta,\eta}$, $\zeta$, $\eta = \pm1$, as defined in Eq.~(\ref{eq:genR})
or~(\ref{eq:genR-exp}), which satisfy the following properties:
\begin{eqnarray}
  (\sigma_{\mu\zeta} \otimes \sigma_{\nu\eta}) \left(\cR^{\zeta,\eta}\right) & =
          & \cR^{\mu,\nu}     \label{eq:genR-cond0} \\
  \tau \circ \Delta^{\eta,\zeta}_{\delta}(A) & = & \cR^{\zeta,\eta}
          \Delta^{\zeta,\eta}_{\delta}(A) \left(\cR^{\zeta,\eta}\right)^{-1}, \qquad 
          A \in \alg,      \label{eq:genR-cond1} \\
  \left(\Delta^{\lambda,\mu}_{\zeta} \otimes \sigma_{\nu\eta}\right) 
          \left(\cR^{\zeta,\eta}\right) & = & \cR^{\lambda,\nu}_{13}\,
          \cR^{\mu,\nu}_{23}, \qquad \left(\sigma_{\lambda\zeta} \otimes
          \Delta^{\mu,\nu}_{\eta}\right) \left(\cR^{\zeta,\eta}\right) = 
          \cR^{\lambda,\nu}_{13}\, \cR^{\lambda,\mu}_{12}.    \label{eq:genR-cond2}  
\end{eqnarray} 
}
\par
%
%
From the results of Proposition~4 or, more simply, by combining
definitions (\ref{eq:extcom-def}), (\ref{eq:genR}), and properties
(\ref{eq:extcom-transf1})--(\ref{eq:extcom-transf3}), (\ref{eq:genR-cond0}), with
Eqs.~(\ref{eq:Rdelta-prop1}), (\ref{eq:Rdelta-prop2}), and~(\ref{eq:Rdelta-prop3}),
we find that the latter can be generalized as follows:
\par
\bigskip
\noindent {\bf Corollary 1}:
{\em The generalized universal \cR-matrix of \alg\ satisfies the relations
\begin{eqnarray}
  \cR^{\zeta,\eta}_{12}\, \cR^{\zeta,\mu}_{13}\, \cR^{\eta,\mu}_{23} & = &
         \cR^{\eta,\mu}_{23}\, \cR^{\zeta,\mu}_{13}\, \cR^{\zeta,\eta}_{12},     
         \label{eq:genR-prop1} \\
  (\epsilon_{\zeta} \otimes \id) \left(\cR^{\zeta,\eta}\right) & = & (\id \otimes
         \epsilon_{\eta}) \left(\cR^{\zeta,\eta}\right) = 1,     \label{eq:genR-prop2} \\
  (S^{\lambda}_{\zeta} \otimes \sigma_{\mu\eta}) \left(\cR^{\zeta,\eta}\right) & =
         & \left(\sigma_{\lambda\zeta} \otimes (S^{\mu}_{\eta})^{-1}\right)
         \left(\cR^{\zeta,\eta}\right) = \left(\cR^{\lambda,\mu}\right)^{-1}.
         \label{eq:genR-prop3}
\end{eqnarray}
}
\par
%
%
In particular, Eq.~(\ref{eq:genR-prop1}) shows that the generalized
\cR-matrix~(\ref{eq:genR-exp}) is a solution of the coloured YBE 
\cite{Murakami}--\cite{Ohtsuki}, where the `colour' parameters $\zeta$,
$\eta$,~$\mu$ take discrete values in the set $\{-1,+1\}$. We therefore propose
to call $\left(\alg,+,m,\iota,\Delta^{\zeta,\eta}_{\delta},\epsilon_{\delta},
S^{\zeta}_{\delta},\cR^{\zeta,\eta};\C\right)$ a two-colour quasitriangular Hopf
algebra over \C.\par
%
%
\section{CONCLUDING REMARKS}   \label{sec:remarks}
In the present paper, we did construct a DQA, denoted as \alg, which has two
series of (N+1)-dimensional unirreps, where $N=0$, 1, 2,~$\ldots$, and we did show
that it can be endowed with a generalized quasitriangular Hopf structure, providing
us with composition laws for all couples of unirreps. This new algebraic structure
was termed a two-colour quasitriangular Hopf algebra because the
corresponding generalized \cR-matrix satisfies the coloured YBE, where the colour
parameters take two discrete values.\par
%
%
It should be noted that various approaches have been previously used to construct
solutions of the coloured YBE \cite{Murakami}--\cite{Kundu}. In the works of Akutsu
and Deguchi~\cite{Akutsu}, and Ge {\it et al}~\cite{Ge}, an infinite-dimensional
representation of sl$_q$(2) was considered and the colour parameter was
introduced as the value of the corresponding Casimir operator. To get
finite-dimensional matrix solutions of the coloured YBE, $q$ had to be restricted to
a root of unity. In the approach pioneered by Burd\'\i k and Hellinger~\cite{Burdik,
Gomez, Basu}, deformations of a non-semisimple Lie algebra, such as gl$_q$(2),
were considered, then the colour parameter was taken as the eigenvalue of the
extra Casimir operator, related with the invariant u(1) subalgebra. Another method,
proposed by Kundu and Basu-Mallick~\cite{Kundu}, used a symmetry transformation
of the YBE (for gl$_q$($N$)) to derive solutions of the coloured one.\par
%
%
Another alternative approach was used in the present work. The colour parameter
now turns out to be related with an involutive automorphism of the algebra
considered. Its two-valuedness is a direct consequence of this property and
contrasts with its continuous character in previous approaches. However, as in the
work of Ge {\it et al}~\cite{Akutsu}, this parameter serves to distinguish between
the representations of the algebra with the same dimension.\par
%
%
It is also worth pointing out that here the colour parameter does not make any
appearance in the algebra defining relations, as it is only needed in the generalized
coalgebraic structure and antipode. This again contrasts with both the coloured
Fadeev-Reshetikhin-Takhtajan algebra and coloured quantum group (generalizing
both GL$_q$(2) and GL$_{p,q}$(2)), recently constructed by
Basu-Mallick~\cite{Basu}, where the colour parameter enters both the algebraic
structure definition and the generator realization, while the coalgebraic
structure remains free from such dependence.\par
%
%
The two-colour quasitriangular Hopf algebra considered here bears some similarity
to the coloured quasitriangular Hopf algebras previously introduced by
Ohtsuki~\cite{Ohtsuki}, which are also characterized by the existence of a coloured
universal \cR-matrix. There are however some differences between both algebraic
structures, the most striking one being the fact that the generalized
comultiplication depends upon two colour parameters in Ref.~\cite{Ohtsuki},
instead of three in the present work.\par
%
%
Construction of other DQA's with a generalized Hopf structure similar to that
considered in the present paper, as well as the investigation of possible
relationships with other coloured algebraic structures, might be some interesting
problems for future study.\par
%
%
\section*{ACKNOWLEDGMENT}
The authors would like to thank the referee for drawing their attention to
Ref.~\cite{Ohtsuki}.\par  
\newpage
%
%
\section*{APPENDIX: DIAGRAMMATIC PROOFS OF EQUATIONS
(\ref{eq:extcom-transf1})--(\ref{eq:extcom-transf3}) AND
(\ref{eq:genHopf1})--(\ref{eq:genHopf3})}
\label{sec:appendix}
\renewcommand{\theequation}{A\arabic{equation}}
\setcounter{equation}{0}
In this Appendix, we prove
Eqs.~(\ref{eq:extcom-transf1})--(\ref{eq:extcom-transf3}) and 
(\ref{eq:genHopf1})--(\ref{eq:genHopf3}) by combining the diagrammatic
definitions~(\ref{eq:diagDelta}) of the generalized coproduct, counit, and antipode
with the action~(\ref{eq:diagsigmaP}) of $\sigma_{\pm}$ on the basic
mapping~$P_{\delta}$, and the diagrammatic representation of standard Hopf
algebra axioms~\cite{Majid}. In the following diagrams, the shorthand notations
\H\ and~\A\ are used for \su\ and \alg, respectively.\par
%
%
Eqs.~(\ref{eq:extcom-transf1})--(\ref{eq:extcom-transf3}) are given by the inner
low rectangular, the outer square, and the outer rectangular diagrams
hereunder, respectively:
\vspace{\abovedisplayskip}
\begin{equation}
  \begin{array}{rcrclcl}
    \H\otimes\H & \stackrel{\Delta}{\la} & \H\,\;\;\; & \stackrel{\ids}{\lra} &
        \,\;\;\;\H & \stackrel{\Delta}{\ra} & \H\otimes\H \\
    {\ss\ids\otimes\ids}\uda\;\;\;\; & & {\ss P_{\delta}}\da\;\;\;\; & & 
        \;\;\;\;\da{\ss P_{\rho}} & & \;\;\;\;\uda{\ss\ids\otimes\ids} \\      
    \H\otimes\H & & \A\,\;\;\; & \stackrel{\sigma_{\rho\delta}}{\ra} & \,\;\;\;\A &
         & \H\otimes\H \\
    {\ss P_{\zeta}\otimes P_{\eta}}\da\;\;\;\; & & {\ss\Delta^{\zeta,\eta}_{\delta}}
         \da\;\;\;\; & & \;\;\;\;{\da}{\ss \Delta^{\mu,\nu}_{\rho}} & &\;\;\;\; 
         \da{\ss P_{\mu}\otimes P_{\nu}} \\
    \A\otimes\A & \stackrel{\ids\otimes\ids}{\lra} & \A\otimes\A &
         \stackrel{\sigma_{\mu\zeta}\otimes\sigma_{\nu\eta}}{\ra} &
         \A\otimes\A & \stackrel{\ids\otimes\ids}{\lra} & \A\otimes\A       
  \end{array}
  \label{eq:A} 
\end{equation}
\vspace{\abovedisplayskip}
\begin{equation}
  \begin{array}{rcrcl}
    \A & \stackrel{\ids}{\lra} & \!\!\A & \stackrel{\ss \epsilon_{\zeta}}{\ra} & \C
        \\
    \ids\uda & & \!\!{\ss P_{\zeta}}\ua & & \uda\ids \\
    \A & & \!\!\H & \stackrel{\ss\epsilon}{\ra} & \C \\
    {\ss\sigma_{\delta\zeta}}\da & & \!\!{\ss P_{\delta}}\da & & \uda\ids
        \\ 
    \A & \stackrel{\ids}{\lra} & \!\!\A & \stackrel{\ss \epsilon_{\delta}}{\ra} & \C
  \end{array}
  \label{eq:B} 
\end{equation}
\vspace{\abovedisplayskip}
\begin{equation}
  \begin{array}{rcrclcl}
    \A & \stackrel{\ids}{\lra} & \A\;\; & \stackrel{\ss S^{\zeta}_{\mu}}{\ra} & \;\;\A
        &\stackrel{\ids}{\lra} & \A \\
    \ids\uda & & {\ss P_{\mu}}\ua\;\; & & \;\;\;\ua{\ss P_{\zeta}} & & \;
        \uda\ids \\
    \A & & \H\;\; & \stackrel{\ss S}\ra & \;\;\H & & \A \\
    {\ss\sigma_{\mu\delta}}\ua & & {\ss P_{\delta}}\da\;\; & & \;\;\;\da {\ss
        P_{\eta}} & & \;\ua{\ss\sigma_{\zeta\eta}} \\
    \A & \stackrel{\ids}{\lra} & \A\;\; & \stackrel{S^{\eta}_{\delta}}\ra & \;\;\A & 
        \stackrel{\ids}{\lra} & \A    
  \end{array}
  \label{eq:C} 
\end{equation}
\vspace{\abovedisplayskip}
\par
%
%
Eq.~(\ref{eq:genHopf1}) corresponds to the outer rectangular diagram:
\vspace{\abovedisplayskip}
\begin{equation}
  \begin{array}{rcrclcl}
    \A & \stackrel{\Delta^{\zeta,\rho}_{\delta}}{\ra} & \;\;\;\;\A\otimes\A & 
         \stackrel{\ids\otimes\ids}{\lra} & \;\;\;\;\A\otimes\A &
         \stackrel{\ids\otimes\ids}{\lra} & \;\;\;\;\A\otimes\A \\
    \ids\uda & & {\ss P_{\zeta}\otimes P_{\rho}}\ua\;\;\;\; & & & &
         \;\;\;\;\;\;\;\;\da{\ss\ids\otimes\Delta^{\eta,\nu}_{\rho}} \\
    \A & & \H\otimes\H & \stackrel{\ids\otimes\Delta}{\ra} & 
         \H\otimes\H\otimes\H & \stackrel{P_{\zeta}\otimes P_{\eta}\otimes
         P_{\nu}}{\ra} & \A\otimes\A\otimes\A \\
    \ids\uda & & {\ss\Delta}\ua\;\;\;\; & & \;\;\;\;\;\;\;\;\ua{\ss\Delta\otimes\ids} 
         & & \;\;\;\;\;\;\;\;\ua{\ss\Delta^{\zeta,\eta}_{\mu}\otimes\ids} \\
    \A & \stackrel{P_{\delta}}{\la} & \H\,\;\;\; & \stackrel{\Delta}{\ra} &
         \;\;\;\;\H\otimes\H & & \;\;\;\;\A\otimes\A \\
    \ids\uda & & & & \;\;\;\;\;\;\;\;\da{\ss P_{\mu}\otimes P_{\nu}} & & 
         \;\;\;\;\;\;\;\;\uda{\ss\ids\otimes\ids} \\
    \A & \stackrel{\ids}{\lra} & \A\,\;\;\; &
         \stackrel{\Delta^{\mu,\nu}_{\delta}}{\ra} & \;\;\;\;\A\otimes\A &
         \stackrel{\ids\otimes\ids}{\lra} & \;\;\;\;\A\otimes\A    
  \end{array}
  \label{eq:D}
\end{equation}
\vspace{\abovedisplayskip}
\par
%
%
The two parts of Eq.~(\ref{eq:genHopf2}) are given by the outer rectangular
diagrams:
\vspace{\abovedisplayskip}
\begin{equation}
  \begin{array}{rcrclcl}
    \A & \stackrel{\ids}{\lra} & \A & \stackrel{\cong}{\ra} & \C\otimes\A &
         \stackrel{\ids}{\lra} & \C\otimes\A \\
    \ids\uda & & & & & & \;\;\;\;\uda\ids \\
    \A & \stackrel{P_{\delta}}{\la} & \H & \stackrel{\cong}{\ra} & \C\otimes\H &
         \stackrel{\ids\otimes P_{\delta}}{\ra} & \C\otimes\A \\
    \ids\uda & & \ids\uda & & \;\;\;\;\ua{\ss\epsilon\otimes\ids} & & 
         \;\;\;\;\ua{\ss\epsilon_{\zeta}\otimes\sigma_{\eta\delta}} \\
    \A & \stackrel{P_{\delta}}{\la} & \H & \stackrel{\Delta}{\ra} & \H\otimes\H &
         \stackrel{P_{\zeta}\otimes P_{\eta}}{\ra} & \A\otimes\A \\
    \ids\uda & & & & & & \;\;\;\;\uda\ids \\
    \A & \stackrel{\ids}{\lra} & \A & \stackrel{\Delta^{\zeta,\eta}_{\delta}}{\ra} &
         \A\otimes\A & \stackrel{\ids}{\lra} & \A\otimes\A
  \end{array}
  \label{eq:Ea}
\end{equation}
\vspace{\abovedisplayskip}
and
\vspace{\abovedisplayskip}
\begin{equation}
  \begin{array}{rcrclcl}
    \A & \stackrel{\ids}{\lra} & \A & \stackrel{\cong}{\ra} & \A\otimes\C &
         \stackrel{\ids}{\lra} & \A\otimes\C \\
    \ids\uda & & & & & & \;\;\;\;\uda\ids \\
    \A & \stackrel{P_{\delta}}{\la} & \H & \stackrel{\cong}{\ra} & \H\otimes\C &
         \stackrel{P_{\delta}\otimes\ids}{\ra} & \A\otimes\C \\
    \ids\uda & & \ids\uda & & \;\;\;\;\ua{\ss\ids\otimes\epsilon} & & 
         \;\;\;\;\ua{\ss\sigma_{\zeta\delta}\otimes\epsilon_{\eta}} \\
    \A & \stackrel{P_{\delta}}{\la} & \H & \stackrel{\Delta}{\ra} & \H\otimes\H &
         \stackrel{P_{\zeta}\otimes P_{\eta}}{\ra} & \A\otimes\A \\
    \ids\uda & & & & & & \;\;\;\;\uda\ids \\
    \A & \stackrel{\ids}{\lra} & \A & \stackrel{\Delta^{\zeta,\eta}_{\delta}}{\ra} &
         \A\otimes\A & \stackrel{\ids}{\lra} & \A\otimes\A
  \end{array}
  \label{eq:Eb}
\end{equation}
\vspace{\abovedisplayskip}
\par
%
%
Finally, the two parts of Eq.~(\ref{eq:genHopf3}) correspond to the right-hand lower
rectangular diagrams:
\vspace{\abovedisplayskip}
\begin{equation}
  \begin{array}{rcrclcl}
    \H\;\;\;\; & \stackrel{\ids}{\lra} & \H\;\;\;\; & \stackrel{P_{\mu}}{\ra} & 
         \;\;\;\;\A & \stackrel{\ids}{\lra} & \;\;\;\;\A \\
    \ids\uda\;\;\;\; & & & & & & \;\;\;\;\uda\ids \\
    \H\;\;\;\; & \stackrel{\iota}{\la} & \C\;\;\;\; & \stackrel{\ids}{\lra} & \;\;\;\;\C 
         & \stackrel{\iota}{\ra} & \;\;\;\;\A \\
    \ids\uda\;\;\;\; & & {\ss\epsilon}\ua\;\;\;\; & &
         \;\;\;\;\ua{\ss\epsilon_{\delta}} & & \;\;\;\;\uda\ids \\
    \H\;\;\;\; & & \H\;\;\;\; & \stackrel{P_{\delta}}{\ra} & \;\;\;\;\A & & \;\;\;\;\A
         \\ 
    {\ss m}\ua\;\;\;\; & & {\ss\Delta}\da\;\;\;\; & &
         \;\;\;\;\da{\ss\Delta^{\zeta,\eta}_{\delta}} & & \;\;\;\;\ua{\ss m} \\  
    \H\otimes\H & & \H\otimes\H & \stackrel{P_{\zeta}\otimes P_{\eta}}{\ra} & 
         \A\otimes\A & & \A\otimes\A \\ 
    \ids\uda\;\;\;\; & & {\ss S\otimes\ids}\da\;\;\;\; & & 
         \;\;\;\;\da{\ss S^{\mu}_{\zeta}\otimes\sigma_{\mu\eta}} & & \;\;\;\;\uda\ids
         \\ 
    \H\otimes\H & \stackrel{\ids}{\lra} & \H\otimes\H &
         \stackrel{P_{\mu}\otimes P_{\mu}}{\ra} & \A\otimes\A &
         \stackrel{\ids}{\lra} & \A\otimes\A 
  \end{array}
  \label{eq:Fa}
\end{equation}
\vspace{\abovedisplayskip}
and
\vspace{\abovedisplayskip}
\begin{equation}
  \begin{array}{rcrclcl}
    \H\;\;\;\; & \stackrel{\ids}{\lra} & \H\;\;\;\; & \stackrel{P_{\mu}}{\ra} & 
         \;\;\;\;\A & \stackrel{\ids}{\lra} & \;\;\;\;\A \\
    \ids\uda\;\;\;\; & & & & & & \;\;\;\;\uda\ids \\
    \H\;\;\;\; & \stackrel{\iota}{\la} & \C\;\;\;\; & \stackrel{\ids}{\lra} & \;\;\;\;\C 
         & \stackrel{\iota}{\ra} & \;\;\;\;\A \\
    \ids\uda\;\;\;\; & & {\ss\epsilon}\ua\;\;\;\; & &
         \;\;\;\;\ua{\ss\epsilon_{\delta}} & & \;\;\;\;\uda\ids \\
    \H\;\;\;\; & & \H\;\;\;\; & \stackrel{P_{\delta}}{\ra} & \;\;\;\;\A & & \;\;\;\;\A
         \\ 
    {\ss m}\ua\;\;\;\; & & {\ss\Delta}\da\;\;\;\; & &
         \;\;\;\;\da{\ss\Delta^{\zeta,\eta}_{\delta}} & & \;\;\;\;\ua{\ss m} \\  
    \H\otimes\H & & \H\otimes\H & \stackrel{P_{\zeta}\otimes P_{\eta}}{\ra} & 
         \A\otimes\A & & \A\otimes\A \\ 
    \ids\uda\;\;\;\; & & {\ss \ids\otimes S}\da\;\;\;\; & & 
         \;\;\;\;\da{\ss\sigma_{\mu\zeta}\otimes S^{\mu}_{\eta}} & & \;\;\;\;\uda\ids
         \\ 
    \H\otimes\H & \stackrel{\ids}{\lra} & \H\otimes\H &
         \stackrel{P_{\mu}\otimes P_{\mu}}{\ra} & \A\otimes\A &
         \stackrel{\ids}{\lra} & \A\otimes\A 
  \end{array}
  \label{eq:Fb}
\end{equation}
\newpage
%
%

\end{document}